\documentstyle{amsppt}
\voffset-10mm
\magnification1200
\pagewidth{130mm}
\pageheight{204mm}
\hfuzz=2.5pt\rightskip=0pt plus1pt
\binoppenalty=10000\relpenalty=10000\relax
\TagsOnRight
\loadbold
\nologo
\addto\tenpoint{\normalbaselineskip=1.2\normalbaselineskip\normalbaselines}
\addto\eightpoint{\normalbaselineskip=1.1\normalbaselineskip\normalbaselines}
\let\le\leqslant
\let\ge\geqslant
\let\epsilon\varepsilon
\let\phi\varphi
\let\wt\widetilde
\redefine\d{\roman d}

\define\Res{\operatorname{Res}}
\topmatter
\title
An elementary proof of Ap\'ery's theorem
\endtitle
\author
Wadim Zudilin \rm(Moscow)
\endauthor
\date
\hbox to100mm{\vbox{\hsize=100mm%
\centerline{E-print \tt math.NT/0202159}
\smallskip
\centerline{17 February 2002}
}}
\enddate
\address
\hbox to70mm{\vbox{\hsize=70mm%
\leftline{Moscow Lomonosov State University}
\leftline{Department of Mechanics and Mathematics}
\leftline{Vorobiovy Gory, Moscow 119899 RUSSIA}
\leftline{{\it URL\/}: \tt http://wain.mi.ras.ru/index.html}
}}
\endaddress
\email
{\tt wadim\@ips.ras.ru}
\endemail
\abstract
We present a new `elementary' proof of the irrationality of~$\zeta(3)$
based on some recent `hypergeometric' ideas of Yu.~Nesterenko, T.~Rivoal,
and K.~Ball, and on Zeilberger's algorithm of creative telescoping.
\endabstract
\endtopmatter
\rightheadtext{An elementary proof of Ap\'ery's theorem}
\leftheadtext{W.~Zudilin}
\document

A question of an arithmetic nature
of the values of Riemann's zeta function
$$
\zeta(s):=\sum_{n=1}^\infty\frac1{n^s}
$$
at odd integral points $s=3,5,7,\dots$
looks like a challenge for Number Theory.
An expected answer `{\it each odd zeta value is transcendental\/}'
is still far from being proved. We dispose of a particular
information on the {\it irrationality\/} of odd zeta values,
namely:
\roster
\item"$\bullet$"
$\zeta(3)$ is irrational (R.~Ap\'ery~\cite{Ap}, 1978);
\item"$\bullet$"
infinitely many of the numbers $\zeta(3),\zeta(5),\zeta(7),\dots$
are irrational (T.~Rivoal \cite{Ri1},~\cite{BR}, 2000);
\item"$\bullet$"
each set $\zeta(s+2),\zeta(s+4),\dots,\zeta(8s-3),\zeta(8s-1)$
with odd $s>1$ contains at least one irrational number
(this author~\cite{Zu1},~\cite{Zu2}, 2001);
\item"$\bullet$"
at least one of the four numbers $\zeta(5),\zeta(7),\zeta(9),\zeta(11)$
is irrational (this author~\cite{Zu3},~\cite{Zu4}, 2001).
\endroster
All these results have a {\it classical\/}
well-poised-hypergeometric origin, and we refer
the reader roused the curiosity of this terminology
to the forthcoming works~\cite{Zu4},~\cite{Zu5}, \cite{RZ}
for details. The aim of this note is to prove Ap\'ery's
famous result by `elementary means'.

\proclaim{Ap\'ery's theorem}
The number $\zeta(3)$ is irrational.
\endproclaim

The idea of the following proof
is due to T.~Rivoal \cite{Ri2}, \cite{Ri3},
who mixed approaches of Yu.~Nesterenko \cite{Ne} and K.~Ball,
and our contribution here is to make a use of Zeilberger's algorithm
of creative telescoping in the most elementary manner.

Our starting point is repetition of \cite{Ne, Section~1}.
For each integer $n=0,1,2,\dots$ define the rational function
$$
R_n(t)
:=\biggl(\frac{(t-1)\dotsb(t-n)}{t(t+1)\dotsb(t+n)}\biggr)^2
$$
and denote by~$D_n$ the least common multiple of the numbers
$1,2,\dots,n$ (and $D_0=1$ for completeness).

\proclaim{Lemma 1 \rm(cf\. \cite{Ne, Lemma~1})}
There holds the equality
$$
F_n:=-\sum_{t=1}^\infty R_n'(t)=u_n\zeta(3)-v_n,
\tag1
$$
where $u_n\in\Bbb Z$, $D_n^3v_n\in\Bbb Z$.
\endproclaim

\demo{Proof}
Taking square of the partial-fraction expansion
$$
\frac{(t-1)\dotsb(t-n)}{t(t+1)\dotsb(t+n)}
=\sum_{k=0}^n\frac{(-1)^{n-k}\binom{n+k}n\binom nk}{t+k}
$$
with a help of the relation
$$
\frac1{t+k}\cdot\frac1{t+l}
=\frac1{l-k}\cdot\biggl(\frac1{t+k}-\frac1{t+l}\biggr)
\qquad\text{for}\quad k\ne l,
$$
we arrive at the formula
$$
R_n(t)=\sum_{k=0}^n\biggl(\frac{A_{2k}^{(n)}}{(t+k)^2}
+\frac{A_{1k}^{(n)}}{t+k}\biggr),
$$
with $A_{jk}=A_{jk}^{(n)}$ satisfying the inclusions
$$
A_{2k}=\binom{n+k}n^2\binom nk^2\in\Bbb Z
\quad\text{and}\quad
D_nA_{1k}\in\Bbb Z,
\qquad k=0,1,\dots,n.
\tag2
$$
Furthermore,
$$
\sum_{k=0}^nA_{1k}=\sum_{k=0}^n\Res_{t=-k}R_n(t)
=-\Res_{t=\infty}R_n(t)=0
$$
since $R_n(t)=O(t^{-2})$ as $t\to\infty$, hence the quantity
$$
\align
F_n
&=\sum_{t=1}^\infty\sum_{k=0}^n
\biggl(\frac{2A_{2k}}{(t+k)^3}+\frac{A_{1k}}{(t+k)^2}\biggr)
=\sum_{k=0}^n\sum_{l=k+1}^\infty
\biggl(\frac{2A_{2k}}{l^3}+\frac{A_{1k}}{l^2}\biggr)
\\
&=2\sum_{k=0}^nA_{2k}\biggl(\sum_{l=1}^\infty-\sum_{l=1}^k\biggr)\frac1{l^3}
+\sum_{k=0}^nA_{1k}\biggl(\sum_{l=1}^\infty-\sum_{l=1}^k\biggr)\frac1{l^2}
\endalign
$$
has the desired form~\thetag{1}, with
$$
u_n=2\sum_{k=0}^nA_{2k},
\qquad
v_n=2\sum_{k=0}^nA_{2k}\sum_{l=1}^k\frac1{l^3}
+\sum_{k=0}^nA_{1k}\sum_{l=1}^k\frac1{l^2}.
\tag3
$$
Finally, using the inclusions~\thetag{2} and
$$
D_n^j\cdot\sum_{l=1}^k\frac1{l^j}\in\Bbb Z
\qquad\text{for}\quad k=0,1,\dots,n, \quad j=2,3,
$$
we deduce that $u_n\in\Bbb Z$ and $D_n^3v_n\in\Bbb Z$ as required.
\enddemo

Since
$$
R_0(t)=\frac1{t^2}, \qquad
R_1(t)=\frac1{t^2}+\frac4{(t+1)^2}-\frac4t+\frac4{t+1},
$$
in accordance with formulae~\thetag{3} we find that
$$
F_0=2\zeta(3) \qquad\text{and}\qquad F_1=10\zeta(3)-12.
\tag4
$$

Now, with a help of Zeilberger's algorithm of creative telescoping
\cite{PWZ, Chapter~6}
we get the rational function $S_n(t):=s_n(t)R_n(t)$, where
$$
s_n(t):=4(2n+1)(-2t^2+t+(2n+1)^2),
\tag5
$$
satisfying the following property.

\proclaim{Lemma 2}
For each $n=1,2,\dots$, there holds the identity
$$
(n+1)^3R_{n+1}(t)-(2n+1)(17n^2+17n+5)R_n(t)+n^3R_{n-1}(t)
=S_n(t+1)-S_n(t).
\tag6
$$
\endproclaim

\demo{`One-line' proof}
Divide both sides of~\thetag{6} by $R_n(t)$ and verify
numerically the identity
$$
\align
&
(n+1)^3\biggl(\frac{t-n-1}{t+n+1}\biggr)^2
-(2n+1)(17n^2+17n+5)+
n^3\biggl(\frac{t+n}{t-n}\biggr)^2
\\ &\qquad
=s_n(t+1)\biggl(\frac{t^2}{(t-n)(t+n+1)}\biggr)^2-s_n(t),
\endalign
$$
where $s_n(t)$ is given in~\thetag{5}.
\enddemo

\proclaim{Lemma 3}
The quantity~\thetag{1} satisfies the difference equation
$$
(n+1)^3u_{n+1}-(2n+1)(17n^2+17n+5)u_n+n^3u_n=0
\tag7
$$
for $n=1,2,\dots$\,.
\endproclaim

\demo{Proof}
Since $R_n'(t)=O(t^{-3})$ and $S_n'(t)=O(t^{-2})$,
differentiating identity~\thetag{6} and summing the result
over $t=1,2,\dots$ we arrive at the equality
$$
(n+1)^3F_{n+1}-(2n+1)(17n^2+17n+5)F_n+n^3F_{n-1}
=S_n'(1).
$$
It remains to note that, for $n\ge1$, both functions
$R_n(t)$ and $S_n(t)=s_n(t)R_n(t)$ have second-order zero
at $t=1$. Thus $S_n'(1)=0$ for $n=1,2,\dots$ and we obtain
the desired recurrence~\thetag{7} for the quantity~\thetag{1}.
\enddemo

Consider another rational function
$$
\wt R_n(t)
:=n!^2(2t+n)\frac{(t-1)\dotsb(t-n)\cdot(t+n+1)\dotsb(t+2n)}
{(t(t+1)\dotsb(t+n))^4}
\tag8
$$
and the corresponding hypergeometric series
$$
\wt F_n:=\sum_{t=1}^\infty\wt R_n(t),
\tag9
$$
proposed by K.~Ball.

\proclaim{Lemma 4 \rm(cf\. \cite{BR, the second proof of Lemma~3})}
For each $n=0,1,2,\dots$, there holds the inequality
$$
0<\wt F_n<20(n+1)^4(\sqrt2-1)^{4n}.
\tag10
$$
\endproclaim

\demo{Proof}
Since $\wt R_n(t)=0$ for $t=1,2,\dots,n$ and $\wt R_n(t)>0$ for $t>n$
we deduce that $\wt F_n>0$.

With a help of elementary inequality
$$
\frac1m\cdot\frac{(m+1)^m}{m^{m-1}}
=\biggl(1+\frac1m\biggr)^m
<e
<\biggl(1+\frac1m\biggr)^{m+1}
=\frac1m\cdot\frac{(m+1)^{m+1}}{m^m}
$$
that yields $(m+1)^m/m^{m-1}<em<(m+1)^{m+1}/m^m$ for $m=1,2,\dots$,
we deduce that
$$
e^{-n}\frac{(m+n)^{m+n-1}}{m^{m-1}}
<m(m+1)\dots(m+n-1)
<e^{-n}\frac{(m+n)^{m+n}}{m^m}.
$$
Therefore, for integers $t\ge n+1$,
$$
\align
\wt R_n(t)\cdot\frac{(t+n)^5}{(2t+n)(t+2n)}
&=n!^2\cdot\frac{(t-1)\dotsb(t-n)\cdot(t+n)\dotsb(t+2n-1)}
{(t(t+1)\dotsb(t+n-1))^4}
\\
&<(n+1)^{2(n+1)}\cdot\frac{t^{5t-4}(t+2n)^{t+2n}}
{(t-n)^{t-n}(t+n)^{5(t+n)-4}}
\endalign
$$
and, as a consequence,
$$
\align
\wt R_n(t)\cdot\frac{t^4(t+n)}{(2t+n)(t+2n)(n+1)^2}
&<(n+1)^{2n}\cdot\frac{t^{5t}(t+2n)^{t+2n}}
{(t-n)^{t-n}(t+n)^{5(t+n)}}
\\
&=\biggl(1+\frac1n\biggr)^{2n}\cdot e^{nf(t/n)}
<e^2\cdot\biggl(\sup_{\tau>1}e^{f(\tau)}\biggr)^n,
\tag11
\endalign
$$
where
$$
f(\tau):=\log\frac{\tau^{5\tau}(\tau+2)^{\tau+2}}
{(\tau-1)^{\tau-1}(\tau+1)^{5(\tau+1)}}.
$$
The unique (real) solution~$\tau_0$ of the equation
$$
f'(\tau)=\log\frac{\tau^5(\tau+2)}{(\tau-1)(\tau+1)^5}=0
$$
in the region $\tau>1$ is the zero of the polynomial
$$
\tau^5(\tau+2)-(\tau-1)(\tau+1)^5
=-\biggl(\tau+\frac12\biggr)\biggl(2\biggl(\tau+\frac12\biggr)^4
-5\biggl(\tau+\frac12\biggr)^2-\frac78\biggr),
$$
hence we can determine it explicitly:
$$
\tau_0=-\frac12+\sqrt{\frac54+\sqrt2}.
$$
Thus,
$$
\align
\sup_{\tau>1}f(\tau)
&=f(\tau_0)=f(\tau_0)-\tau_0f'(\tau_0)
=2\log(\tau_0+2)+\log(\tau_0-1)-5\log(\tau_0+1)
\\
&=4\log(\sqrt2-1)
\endalign
$$
and we can continue the estimate~\thetag{11} as follows:
$$
\wt R_n(t)\cdot\frac{t^4(t+n)}{(2t+n)(t+2n)}
<e^2(n+1)^2(\sqrt2-1)^{4n},
\tag12
$$
Finally, we apply the inequality~\thetag{12} to deduce
the required estimate~\thetag{10}:
$$
\align
\wt F_n
&=\sum_{t=n+1}^\infty\wt R_n(t)
<e^2(n+1)^2(\sqrt2-1)^{4n}\sum_{t=n+1}^\infty\frac{(2t+n)(t+2n)}{t^4(t+n)}
\\
&<e^2(n+1)^2(\sqrt2-1)^{4n}\sum_{t=n+1}^\infty
\biggl(\frac2{t^5}+\frac{5n}{t^4}+\frac{2n^2}{t^3}\biggr)
\\
&\le e^2(n+1)^2\bigl(2\zeta(5)+5n\zeta(4)+2n^2\zeta(3)\bigr)
(\sqrt2-1)^{4n}
<20(n+1)^4(\sqrt2-1)^{4n}.
\endalign
$$
This completes the proof.
\enddemo

For the rational function~\thetag{8} we obtain Zeilberger's certificate
$$
\align
\wt S_n(t)
&:=\frac{\wt R_n(t)}{(2t+n)(t+2n-1)(t+2n)}\cdot
\bigl(-t^6-(8n-1)t^5+(4n^2+27n+5)t^4
\\ &\qquad
+2n(67n^2+71n+15)t^3
+(358n^4+339n^3+76n^2-7n-3)t^2
\\ &\qquad
+(384n^5+396n^4+97n^3-29n^2-17n-2)t
\\ &\qquad
+n(153n^5+183n^4+50n^3-30n^2-22n-4)\bigr).
\tag13
\endalign
$$

\proclaim{Lemma 5}
For each $n=1,2,\dots$, there holds the identity
$$
(n+1)^3\wt R_{n+1}(t)-(2n+1)(17n^2+17n+5)\wt R_n(t)+n^3\wt R_{n-1}(t)
=\wt S_n(t+1)-\wt S_n(t).
\tag14
$$
\endproclaim

\demo{`One-line' proof}
Divide both sides of~\thetag{14} by $\wt R_n(t)$ and verify
the reduced identity.
\enddemo

\proclaim{Lemma 6}
The quantity~\thetag{9} satisfies the difference equation~\thetag{7}
for $n=1,2,\dots$\,.
\endproclaim

\demo{Proof}
Since $\wt R_n(t)=O(t^{-5})$ and $\wt S_n(t)=O(t^{-2})$
as $t\to\infty$ for $n\ge1$, summation of equalities~\thetag{14}
over $t=1,2,\dots$ yields the relation
$$
(n+1)^3\wt F_{n+1}-(2n+1)(17n^2+17n+5)\wt F_n+n^3\wt F_{n-1}
=-\wt S_n(1).
$$
It remains to note that, for $n\ge1$, both functions~\thetag{8}
and~\thetag{13} have zero
at $t=1$. Thus $\wt S_n(1)=0$ for $n=1,2,\dots$ and we obtain
the desired recurrence~\thetag{7} for the quantity~\thetag{9}.
\enddemo

\proclaim{Lemma 7}
For each $n=0,1,2,\dots$, the quantities~\thetag{1} and~\thetag{9}
coincide.
\endproclaim

\demo{Proof}
Since both $F_n$ and $\wt F_n$ satisfy the same second-order
difference equation~\thetag{7}, we have to verify that
$F_0=\wt F_0$ and $F_1=\wt F_1$. Direct calculations show that
$$
\wt R_0(t)=\frac2{t^3},
\qquad
\wt R_1(t)
=-\frac2{t^4}+\frac2{(t+1)^4}
+\frac5{t^3}+\frac5{(t+1)^3}
-\frac5{t^2}+\frac5{(t+1)^2},
$$
hence $\wt F_0=2\zeta(3)$ and $\wt F_1=10\zeta(3)-12$, and
comparison of this result with~\thetag{4} yields
the desired coincidence.
\enddemo

\demo{Proof of Ap\'ery's theorem}
Suppose, on the contrary, that $\zeta(3)=p/q$, where $p$ and~$q$
are positive integers. Then, using a trivial bound
$D_n<3^n$, we deduce that, for each $n=0,1,2,\dots$,
the integer $qD_n^3F_n=D_n^3u_np-D_n^3v_nq$ satisfies the estimate
$$
0<qD_n^3F_n<20q(n+1)^43^{3n}(\sqrt2-1)^{4n}
\tag15
$$
that is not possible since $3^3(\sqrt2-1)^4=0.7948\hdots<1$ and
the right-hand side of~\thetag{15} is less than~$1$ for a
sufficiently large integer~$n$.
This contradiction completes the proof of the theorem.
\enddemo

Inspite of its elementary arguments, our proof of Ap\'ery's
theorem does not look simpler than the original
(also elementary) Ap\'ery's proof
well-explained in A.~van der Poorten's informal report~\cite{Po},
or (almost elementary) Beukers's proof~\cite{Be} by means
of Legendre polynomials and multiple integrals. We want to mention
that our way to deduce the recursion~\thetag{7} for the sequence
$F_n$ as well as for the coefficients $u_n,v_n$%
\footnote"${}^\ddag$"{Hint:
multiply both sides of~\thetag{6} by $(t+k)^2$,
substitute $t=-k$ and sum over {\it all\/} integers~$k$ to show
that the sequence $u_n$ satisfies the difference equation~\thetag{7};
then $v_n=u_n\zeta(3)-F_n$ also satisfies it.}
slightly differs from those considered in~\cite{Po, Section~8}
and~\cite{Ze, Section~13} although it is based on the same
algorithm of creative telescoping. This algorithm and the above
scheme allow us~\cite{Zu5},~\cite{Zu6} to obtain
Ap\'ery-like difference equations
for $\zeta(4)$ and Calalan's constant.

The fact that $\wt F_n=\wt u_n\zeta(3)-\wt v_n$ with
$D_n\wt u_n,D_n^4\wt v_n\in\Bbb Z$ was first discovered
by K.~Ball; the proof follows lines of the proof of Lemma~1
and vanishing the coefficients for $\zeta(4)$ and $\zeta(2)$
is due to a well-poised origin of the series~\thetag{9}.
An open question of T.~Rivoal
here is to get the better inclusions $\wt u_n,D_n^3\wt v_n\in\Bbb Z$
by elementary means without going back to Ap\'ery's series~\thetag{1}.
A solution of this question accompanied with Ball's Lemma~4
can bring the `most elementary' proof of Ap\'ery's theorem.

Lemma~7 can be proved by specialization
of Bailey's identity~\cite{Ba, Section~6.3, formula~(2)}
$$
\align
&
{}_7\!F_6\biggl(\matrix\format&\,\c\\
a, & 1+\frac12a, &     b, &     c, &     d, &     e, &     f \\
   &   \frac12a, & 1+a-b, & 1+a-c, & 1+a-d, & 1+a-e, & 1+a-f
\endmatrix\biggm|1\biggr)
\\ &\qquad
=\frac{\Gamma(1+a-b)\,\Gamma(1+a-c)\,\Gamma(1+a-d)\,
\Gamma(1+a-e)\,\Gamma(1+a-f)}
{\aligned
\Gamma(1+a)\,\Gamma(b)\,\Gamma(c)\,\Gamma(d)\,
\Gamma(1+a-b-c)\,\Gamma(1+a-b-d)\,
\qquad\quad \\ \vspace{-.3\baselineskip} \times
\Gamma(1+a-c-d)\,\Gamma(1+a-e-f)
\endaligned}
\\ &\qquad\quad\times
\frac1{2\pi i}\int_{-i\infty}^{i\infty}
\frac{\aligned
\Gamma(b+t)\,\Gamma(c+t)\,\Gamma(d+t)\,\Gamma(1+a-e-f+t)\,
\qquad\quad \\ \vspace{-.3\baselineskip} \times
\Gamma(1+a-b-c-d-t)\,\Gamma(-t)
\endaligned}
{\Gamma(1+a-e+t)\,\Gamma(1+a-f+t)}
\,\d t,
\tag16
\endalign
$$
provided that the very-well-poised hypergeometric series
on the left-hand side converges.
Namely, taking $a=3n+2$ and $b=c=d=e=f=n+1$ in~\thetag{16}
we obtain Ball's sequence~\thetag{9}
on the left and Ap\'ery's sequence~\thetag{1} on the right
(for the last fact see~\cite{Ne, Lemma~2}). Identity~\thetag{16}
can be put forward for an explanation how the permutation
group from~\cite{RV} for linear forms in~$1$ and~$\zeta(3)$
appears (see \cite{Zu5, Sections~4 and~5 for details}).



\Refs
\widestnumber\key{WW}

\ref\key Ap
\by R.~Ap\'ery
\paper Irrationalit\'e de $\zeta(2)$ et $\zeta(3)$
\jour Ast\'erisque
\vol61
\yr1979
\pages11--13
\endref

\ref\key Ba
\by W.\,N.~Bailey
\book Generalized hypergeometric series
\bookinfo Cambridge Math. Tracts
\vol32
\publ Cambridge Univ. Press
\publaddr Cambridge
\yr1935
\moreref
\bookinfo 2nd reprinted edition
\publaddr New York--London
\publ Stechert-Hafner
\yr1964
\endref

\ref\key BR
\by K.~Ball and T.~Rivoal
\paper Irrationalit\'e d'une infinit\'e de valeurs
de la fonction z\^eta aux entiers impairs
\jour Invent. Math.
\vol146
\yr2001
\issue1
\pages193--207
\endref

\ref\key Be
\by F.~Beukers
\paper A note on the irrationality of~$\zeta(2)$ and~$\zeta(3)$
\jour Bull. London Math. Soc.
\vol11
\issue3
\yr1979
\pages268--272
\endref

\ref\key Ne
\by Yu.\,V.~Nesterenko
\paper A few remarks on~$\zeta(3)$
\jour Mat. Zametki [Math. Notes]
\vol59
\yr1996
\issue6
\pages865--880
\endref

\ref\key PWZ
\by M.~Petkov\v sek, H.\,S.~Wilf, and D.~Zeilberger
\book $A=B$
\publaddr Wellesley, M.A.
\publ A.\,K.~Peters, Ltd.
\yr1997
\endref

\ref\key Po
\by A.~van der Poorten
\paper A proof that Euler missed...
Ap\'ery's proof of the irrationality of~$\zeta(3)$
\paperinfo An informal report
\jour Math. Intelligencer
\vol1
\issue4
\yr1978/79
\pages195--203
\endref

\ref\key RV
\by G.~Rhin and C.~Viola
\paper The group structure for~$\zeta(3)$
\jour Acta Arith.
\vol97
\issue3
\yr2001
\pages269--293
\endref

\ref\key Ri1
\by T.~Rivoal
\paper La fonction z\^eta de Riemann prend une infinit\'e
de valeurs irrationnelles aux entiers impairs
\jour C.~R. Acad. Sci. Paris S\'er.~I Math.
\vol331
\yr2000
\issue4
\pages267--270
\moreref
\inbook E-print {\tt math.NT/0008051}
\endref

\ref\key Ri2
\by T.~Rivoal
\book Propri\'et\'es diophantinnes des valeurs
de la fonction z\^eta de Riemann aux entiers impairs
\bookinfo Th\`ese de Doctorat
\publ Univ. de Caen
\publaddr Caen
\yr2001
\endref

\ref\key Ri3
\by T.~Rivoal
\paper S\'eries hyperg\'eom\'etriques et irrationalit\'e
des valeurs de la fonction z\^eta
\inbook Journ\'ees arithm\'etiques (Lille, July, 2001)
\yr2002
\toappear
\endref

\ref\key RZ
\by T.~Rivoal and W.~Zudilin
\paper Diophantine properties of numbers related to Catalan's constant
\inbook Pr\'epublication de l'Institut de Math. de Jussieu,
no.~315 (Janvier 2002)
\finalinfo submitted for publication
\endref

\ref\key Ze
\by D.~Zeilberger
\paper Closed form (pun intended!)
\paperinfo A tribute to Emil Crosswald
``Number theory and related analysis''
\eds M.~Knopp and M.~Sheingorn
\inbook Contemporary Math.
\vol143
\publaddr Providence, R.I.
\publ Amer. Math. Soc.
\yr1993
\pages579--607
\endref

\ref\key Zu1
\by W.~Zudilin
\paper Irrationality of values of zeta-function
\inbook Contemporary Research in Mathematics and Mechanics
\bookinfo Proceedings of the XXIII Conference of Young Scientists
of the Department of Mechanics and Mathematics
(Moscow State University, April~9--14, 2001)
\publaddr Moscow
\publ Publ. Dept. Mech. Math. MSU
\yr2001, Part~2
\pages 127--135
\moreref
\inbook English transl.,
E-print {\tt math.NT/\allowlinebreak0104249}
\endref

\ref\key Zu2
\by W.~Zudilin
\paper Irrationality of values of the Riemann zeta function
\jour Izv. Ross. Akad. Nauk Ser. Mat.
[Russian Acad. Sci. Izv. Math.]
\vol66
\yr2002
\issue3
\endref

\ref\key Zu3
\by W.\,V.~Zudilin
\paper One of the numbers
$\zeta(5),\zeta(7),\zeta(9),\zeta(11)$ is irrational
\jour Uspekhi Mat. Nauk [Russian Math. Surveys]
\vol56
\yr2001
\issue4
\pages149--150
\endref

\ref\key Zu4
\by W.~Zudilin
\paper Arithmetic of linear forms involving odd zeta values
\jour Preprint (August 2001)
\finalinfo submitted for publication
\endref

\ref\key Zu5
\by W.~Zudilin
\paper Difference equation and permutation group for~$\zeta(4)$
\paperinfo Actes des 12\`emes rencontres arithm\'etiques de Caen
(June 29--30, 2001)
\jour J. Th\'eorie Nombres Bordeaux
\yr2002
\toappear
\endref

\ref\key Zu6
\by W.~Zudilin
\paper Ap\'ery-like difference equation for Catalan's constant
\jour Preprint (January 2002),
E-print {\tt math.NT/\allowlinebreak0201024}
\finalinfo submitted for publication
\endref

\endRefs
\enddocument